 \documentclass[final,1p,times]{elsarticle}

\usepackage{amssymb}

\usepackage[figuresright]{rotating}
\usepackage{amsmath}
\usepackage{amscd}
\usepackage{wrapfig}
\usepackage{blkarray}
\usepackage{hyperref}
\usepackage{color}
\usepackage{kbordermatrix}

\usepackage{amsfonts}
\usepackage{mathrsfs}
\usepackage{stmaryrd}

\usepackage{epsfig}
\usepackage{epstopdf}
\usepackage{graphicx}
\usepackage{caption}
\usepackage{subcaption}
\usepackage{float}

\usepackage{algorithm}
\usepackage{algcompatible}
\usepackage{multirow}

\usepackage{dsfont}
\usepackage{url}

\newtheorem{example}{Example}[section]

\newtheorem{definition}{Definition}[section]
\newtheorem{theorem}{Theorem}[section]

\newtheorem{lemma}[theorem]{Lemma}
\newtheorem{remark}[theorem]{Remark}

\newlength{\ColWidth}

\newcommand{\arrayvline}{\hspace*{\arraycolsep}\vline\hspace*{-\arraycolsep}}

\begin{document}

\begin{frontmatter}



\title{Efficient algorithm for computing large scale systems of differential algebraic equations}


%

\author[Casit,LiU]{Xiaolin Qin}
\ead{qinxl@casit.ac.cn}
\author[Casit]{Juan Tang \corref{cor1}}
\ead{tangjuan0822@gmail.com}
\author[Casit]{Yong Feng}
\author[Buas]{Bernhard Bachmann}
\author[LiU]{Peter Fritzson}

\cortext[cor1]{Corresponding author}
\address[Casit]{Chengdu Institute of Computer Applications, Chinese Academy of Sciences, Chengdu 610041, China}
\address[LiU]{Department of Computer and Information Science, Link\"{o}ping University, Link\"{o}ping SE-581 83, Sweden}
\address[Buas]{Department of Mathematics and Engineering, Bielefeld University of Applied Sciences, Bielefeld D-33609, Germany}

\begin{abstract}
In many mathematical models of physical phenomenons and engineering fields, such as electrical circuits or mechanical multibody systems, which generate the differential algebraic equations (DAEs) systems naturally. In general, the feature of DAEs is a sparse large scale system of fully nonlinear and high index. To make use of its sparsity, this paper provides a simple and efficient algorithm for computing the large scale DAEs system. We exploit the shortest augmenting path algorithm for finding maximum value transversal (MVT) as well as block triangular forms (BTF). We also present the extended signature matrix method with the block fixed point iteration and its complexity results. Furthermore, a range of nontrivial problems are demonstrated by our algorithm.

\end{abstract}

\begin{keyword}
Differential algebraic equations, sparsity, shortest augmenting path, block triangular forms, structural analysis
\end{keyword}

\end{frontmatter}


\section{Introduction}
The problem of differential algebraic equations (DAEs) system solving is fundamental in modelling many equation-based models of physical
phenomenons and engineering fields, such as electric circuits
\cite{jour15,jour16}, mechanical systems \cite{book3}, spacecraft dynamics \cite{jour17}, chemical engineering \cite{VB2000}, and many other areas. Generally, DAEs can be produced very large scale system of fully nonlinear and higher index in practice. However, most of the algorithms treat the low index case or consider solutions of linear systems \cite{KM2006, M1992, Pantelides1988, TI2008}. The index is a notion used in the theory of DAEs for measuring the distance from a DAE to its related ordinary differential equation (ODE). It is well known that it is direct numerical computations difficult to solve a high index DAE. In particular, it may only solve some special classes of DAEs by the direct numerical solution \cite{EH2009, L2014}. 

Index reduction techniques can be used to remedy the drawback of direct numerical computation \cite{BCP1996}. Pantelides' method \cite{Pantelides1988} gives a systematic way to reduce high index systems of DAEs to lower index one, by selectively adding differentiated forms of the equations already present in the system. In \cite{DZC2008}, Ding et al. developed the weighted bipartite algorithm based on the minimally structurally singular subset, which is similar to the Pantelides' method. However, the algorithms can succeed yet not correctly in some instances \cite{Reissig2000}. Campbell's derivative array \cite{C1993} needs to be computationally expensive especially for computing the singular value decomposition of the Jacobian of the derivative array equations using nonlinear singular least squares methods. Signature matrix method (also called $\Sigma$-method) \cite{Pryce2001} is based on solving an assignment problem, which can be formulated as an integer linear programming problem. In \cite{Pryce2001}, Pryce proved that $\Sigma$-method is equivalent to the famous method of Pantelides' algorithm, and in particular computes the same structural index. However, the nice feature of $\Sigma$ method is a simple and straightforward method for analyzing the structure of DAEs of any order, not just first order.

In particular, large scale and high index DAEs with fully nonlinear systems are now routine that such models are built using
interactive design systems based on the Modelica language \cite{CG2011, F2015}. In addition, the sparsity pattern of DAEs can arise in most actual applications \cite{Pryce2001, QWFR2013}. In \cite{FKF2012}, Frenkel et al. gave a survey on appropriate matching algorithms based on the augmenting paths and push-relabel algorithm by translating Modelica models for large scale systems of DAEs. More recently, Wu et al. \cite{WRI2009} generalized the $\Sigma$-method to the square and $t$-dominated partial differential equations (PDEs) systems. Pryce et al. \cite{Pryce2014} generalized the $\Sigma$-method for constructing a block triangular form (BTF) of the DAEs and exploiting to solve it efficiently in a block-wise manner. In \cite{TWQF2014}, Tang et al. proposed the block fixed-point iteration with parameter method for DAEs based on its block upper triangular structure. However, the essential task is to solve the linear assignment problem for finding a maximum value transversal (MVT), which is a large part of the cost for index reduction of DAEs solving. Pryce et al. mentioned only in their work using Cao's \emph{Matlab} implementation \cite{C2013} of the shortest augmenting path algorithm of Jonker and Volgenant in \cite{JV1987}. We focus on solving in the large scale and high index cases in order to provide the shortest augmenting path algorithm for finding an MVT and an extended signature matrix method. The problem is also closely related to computing the block triangular form of a sparse matrix and linear assignment problems over integer.

Our approach is based on signature matrix method and modified Dijkstra's shortest path method. Our fundamental tool is the block triangularization of a sparse matrix; we exploit recent advances in linear assignment problem solving, which is equivalent to finding a maximum weight perfect matching in a bipartite graph of signature matrix in $\Sigma$-method, and we adapt the block fixed-point iteration with parameter for the canonical offsets techniques. Currently, we are working on the theoretical foundation and implementation of these methods on \emph{Maple} platform. Another direction for future work is to exploit the fact that our algorithms are expressed in OpenModelica solvers.

The rest of this paper is organized as follows. The next section introduces our purpose and the shortest augmenting paths based algorithms, and presents an improved algorithm for the block triangularization for DAEs system. Section 3 describes the extended signature matrix method for the structural analysis of large scale DAEs system and gives its complexity results. The following section shows our algorithm for an actual industrial example and some experimental results. The final section makes conclusions.

\section{Preliminary results }

\subsection{Purpose} \label{sec:purp}
We consider an input DAEs system in $n$ dependent variables $x_{j}=x_{j}(t)$
with $t$ a scalar independent variable, of the general form
\begin{equation}
f_{i}(t,the \ x_j\ and\ derivatives\ of\ them)=0,\ 1 \leq i, j \leq n.
\end{equation}
The $f_i$ are assumed suitably smooth, and the derivatives of $x_j$ are arbitrary order. In general, signature matrix method is an effective preprocessing algorithm for the small and middle scale DAEs system. First, it needs to form the $n \times n$ signature matrix $\Sigma =
(\sigma_{ij})$ of the DAEs, where
\begin{center} 
${\sigma_{ij}}=\begin{cases}
 highest \ order \ of \ derivative \ to \ which \ the \ variable \ x_j \ appears \ in \ equation \ f_i, \\
or \ -\infty \ if \ the \ variable \ does \ not \ occur.
\end{cases}$
\end{center}
Then, taken the analysis procedure of $\Sigma$ as a linear assignment problem is to seek the offsets of the DAEs, that is, the number of differentiations of $f_i$. It can be formulated by the following primal problem:
\begin{eqnarray}\label{Lap:Primal}
\begin{array}{ll}
  Maximize & z =\sum\limits_{(i,j)\in S} \sigma_{ij}\xi_{ij}, \label{equ:a}   \\
  subject \ to & \sum\limits_{\{j|(i,j)\in S\}} \xi_{ij}=1 \ \ \ \ \forall\ i =1, 2, \cdots, n, \label{equ:b}\\
  & \sum\limits_{\{i|(i,j)\in S\}} \xi_{ij}=1 \ \ \ \ \forall\ j = 1, 2, \cdots, n, \label{equ:c}\\
  & \xi_{ij} \in \{0, 1 \} \ \ \ \ \ \ \ \ \ \ \forall\ (i,j) \in S. \label{equ:d}
\end{array}
\end{eqnarray}

Note that the state variable $\xi_{ij}$ only be defined over the sparsity pattern of the problem
\begin{equation}
  S=sparse(\Sigma)=\{(i,j)| \sigma_{ij} > -\infty \}.
\end{equation}
It can be also defined on an undirected bipartite graph, in which case an assignment is a perfect matching. Given a bipartite graph $G(\Sigma) = (F, X, e)$, where $F$ is the set vertices corresponding to the rows of $\Sigma$, $X$ is the set vertices corresponding to the columns of $\Sigma$, and $e$ is the set of edges corresponding to the non-negatively infinite in $\Sigma$, $|F| = |X| = n$, $|\cdot|$ denotes the cardinality of a set. In this paper, our goal is to handle large scale systems with $n$ involving thousands and even more.

\subsection{Block triangularization for DAEs system}
\label{sec:btf}
As we have encountered with the increasingly large problems, an important preprocessing manipulation is the block triangularization of the system,
which allows to split the overall system into subsystems which can be solved in the sequence or parallelization \cite{CM1998}. Pothen et al. \cite{PF1990} described implementations of algorithms to compute the block triangular form by permuting the rows and columns of a rectangular or square, unsymmetric sparse matrix. It is equivalent to computing a canonical decomposition of bipartite graphs known as the Dulmage-Mendelsohn decomposition. Considering the index reduction for DAEs system, the block triangularization algorithm can be directly performed on the signature matrix $\Sigma$. The block triangular form of $\Sigma$ can be generated by the block triangularization of the incidence matrix of the sparse pattern $S$ for a given DAEs system. The incidence matrix is defined as follows:
\begin{definition}
Let $A =[a_{ij}]$ be an incidence matrix of the sparse pattern $S$ for a given DAEs system, whose rows and columns represent the $F$ and $X$ as above, and element $a_{ij}$ is $1$ if $(i,j) \in S$ and is $0$ otherwise.
\end{definition}
The associated bipartite graph $G(A)$ for the incidence matrix $A$ is defined as follows:
\begin{definition}
Let $G(A) = (F, X, E)$ be the associated undirected bipartite graph for the
incidence matrix $A$ of a given DAEs system, where $F$ and $X$ are defined as above, one for each equation $f_i$, and other for each variable $x_j$ respectively; and $(f_{i},x_{j})$ belongs to $E$ if and only if $a_{ij}{\neq}0$.
\end{definition}

\begin{definition}
A subset $M$ of $E$ is called a matching if a vertex in $G(A)$ is incident to at most one edge in $M$. A matching $M$ is called
maximal, if no other matching $M'\supset  M$ exists, which is called maximum $if |M| \geq |M'|$ for every matching $M'$. Furthermore, if $|M| = |F| = |X|$, $M$ is called a perfect matching.
\end{definition}

\begin{lemma}(\cite{LP2009})\label{lem:hall}
A bipartite graph $G(A)$ with vertex sets $F$ and $X$ contains a perfect matching from $F$ to $X$ if and only if it satisfies Hall's condition
\begin{eqnarray*}
|\Gamma(\textbf{f})|\geq |\textbf{f}| \ for \ every \ \textbf{f} \subset F,
\end{eqnarray*}
where $\Gamma(\textbf{f})=\{ x_j \in X, (f_i, x_j)\in \ E \ for \ some \ f_i \in F \} \subset X$.
\end{lemma}

A walk is a sequence of nodes $v_0, v_1, \cdots, v_{n-1} \in F \bigcup X$
such that $(v_i, v_{i+1}) \in E$ for $i =0, 1, \cdots, n-2$.
Noted that edges or nodes can be repeated in a walk. An alternating walk is a walk with alternate edges in a matching $M$.
An alternating path is an alternating walk with no repeated nodes. With respect to $M$, we can define the following sets:
\begin{equation*}
\begin{array}{l}
VF=\{$F-nodes$ \ reachable \ by \ alternating \ path \ from \ some \ unmatched\ $F-node$\} \\
HF=\{$F-nodes$ \ reachable\ by\ alternating\ path\ from\ some\ unmatched\ $X-node$ \} \\
SF=F\setminus(VF \bigcup HF) \\
VX=\{$X-nodes$\ reachable\ by\ alternating\ path\ from\ some\ unmatched\ $F-node$\} \\
HX=\{$X-nodes$\ reachable\ by\ alternating\ path\ from\ some\ unmatched\ $X-node$\} \\
SX=X\setminus(VX \bigcup HX)
\end{array}
\end{equation*}

If there exists a perfect matching for $G(A)$, the directed graph $G_{d}(A)$
of $G(A)$ is made by the rule: let non-matching edges in $G(A)$ be directed from $X$-nodes to $F$-nodes,
matching edges shrunk into single nodes, and the nodes identified with $F$. In order to demonstrate the benefits of parallel processing for index reduction, we present an improved algorithm for the block triangularization of the incidence matrix $A$. The main steps are as follows:
\begin{description}
  \item[step 1] Find the diagonal sub-blocks of $A$ from the connected components of $G(A)$. Assume it generates
  $p$ diagonal submatrices ${A_{1}, A_{2}, ... ,A_{p}}$. It is helpful to use parallel computing for each $A_i$.
  \begin{equation*}
\begin{array}{c}
A\stackrel{Step1}{\longrightarrow}
\begin{bmatrix}
A_1 \\
& A_2 & & \text{{\huge{0}}}\\
& & \ddots \\
& \text{{\huge{0}}} & & A_{p-1}\\
& & & & A_{p}
\end{bmatrix}
\end{array}
\end{equation*}

  \item[step 2] If $A$ has a perfect matching, then it only has the square $A_s$, otherwise the underdetermined $A_h$ and overdetermined $A_v$ will be present. For each $A_{i}$, find a maximum matching ${M}_{i}$ in the associated graph $G(A_{i})$.

  \item[step 3] With respect to $M_{i}$, partition $F_{i}$ into the sets $VF_{i}$, $SF_{i}$, $HF_{i}$ using the rule of Dulmage-Mendelsohn decomposition; partition $X_{i}$ into the sets $VX_{i}$, $SX_{i}$, $HX_{i}$ similarly, where $F_i$ and $X_i$ are the block form of $F$ and $X$ respectively, $'\times'$ denotes a possibly nonzero matrix of appropriate dimensions.
      \begin{equation*}
\begin{array}{c}
A_i \stackrel{Step3}{\longrightarrow}
\begin{blockarray}{cccc}
~    & HX_i       & SX_i       & VX_i \\
\begin{block}{c[ccc]}
  HF_i & A_{i_h}    & \times          & \times    \\
 SF_i & \textbf{0} & A_{i_s}    & \times    \\
  VF_i & \textbf{0} & \textbf{0} & A_{i_v} \\
\end{block}
\end{blockarray}
\end{array}
\end{equation*}
  \item[step 4] Using Tarjan's depth-first search algorithm, find the block upper triangular form of the square submatrix $ A_{i_{s}}$ ($A_{i_{v}}$ and $A_{i_{h}}$ are not square) by finding strong components in the associated directed graph $ G_{d}(A_{i_{s}})$.
       Assume it produces $A_{i_{s},1}, A_{i_{s},2}, ..., A_{i_{s},n_{i}}$.
       \begin{equation*}
\begin{array}{c}
A_{i_s} \stackrel{Step4}{\longrightarrow}
\begin{blockarray}{ccccc}
     ~ & SX_{i,1}  & SX_{i,2}  & \cdots   & SX_{i,n_i}\\
\begin{block}{c[cccc]}
SF_{i,1} & A_{i_s,1} & \times         & \cdots    & \times          \\
SF_{i,2} &\textbf{0} & A_{i_s,2} & \cdots    & \times         \\
 \vdots & \vdots    & \vdots    & \ddots    & \vdots     \\
SF_{i,n_i}&\textbf{0}& \textbf{0}& \cdots    & A_{i_s,n_i}  \\
\end{block}
\end{blockarray}
\end{array}
\end{equation*}
\end{description}
\begin{remark}
Some problems may only require one step of the computation in step 1, one at a time; others can take advantage of the further decomposition. In addition, if ${M}_{i}$ is perfect matching($A_i$ contains at least one
transversal), then $HX_i$,$VX_i$,$HF_i$ and $VF_i$ are empty, that is,
the step $3$ can be skipped and $A_i=A_{i_s}$. Furthermore, our algorithm can be generalized to the underdetermined and overdetermined signature matrix of sparse pattern $S$ case.
\end{remark}
\begin{remark}\label{lem:irreducibility}
\begin{description}
  \item[(a)] The following are equivalent:
  \begin{description}
    \item[(i)] $A$ is structurally nonsingular (contains at least one transversal).
    \item[(ii)] $A$ has the Hall's condition from Lemma \ref{lem:hall}.
 \end{description}
  \item[(b)] The following are equivalent:
  \begin{description}
    \item[(i)] $A$ is irreducible.
    \item[(ii)] $A$ has the strong Hall's condition, which satisfies for all proper subsets $|\Gamma(\textbf{f})|\geq |\textbf{f}|+1$.
    \item[(iii)] Every element of $A$ is on a transversal, and the bipartite graph of $A$ is connected.
  \end{description}
\end{description}
\end{remark}
\begin{theorem}
The above algorithm works correctly as specified and its
complexity as follows:\\
(a) If $A$ is irreducible, then it only performances in step 1, that is, $\textit{O}(n + \tau)$ operations,\\
(b) If $A$ is reducible, then it needs the whole steps, that is, $\textit{O}(pN\Phi)$ operations,\\
where $N$ and $\Phi$ are the constant, as defined below.
\end{theorem}
\begin{proof}
Correctness of the algorithm follows from the basic idea. The essential tasks of above algorithm are to permute the given matrix to block triangular form. A permuted block form of the DAEs system is generated from forming permutations $\tilde{F}$ and $\tilde{X}$ of the equations and variables, which denote the block form $(F_1, F_2, \cdots, F_p)$ and $(X_1, X_2, \cdots, X_p)$ with $n_1 + n_2 + \cdots + n_p = n$, where $|F_i| = |X_i| =n_i, i=1, 2, \cdots, p$.

From the description of algorithm, we observe that there are four major steps on time complexity. In step 1, we can use a version of recursive depth-first search to find the connected components in the run time $\textit{O}(n + \tau)$, where $\tau$ is the number of nonzero entries in $A$ ($\tau = \tau_1 + \tau_2 + \cdots + \tau_p $). Suppose $n_i = N$, $\tau_i =\Phi$  \footnote{$N, \Phi$ are defined by the same way for the rest of this paper.}, for each $i$, that it, $n = p \cdot N$, $\tau = p \cdot \Phi$, a $\textit{O}(N\Phi)$ algorithm for each $A_i (i=1, 2,\cdots, p)$ is considered in step 2, which can be parallel computing. It terminates with a maximum matching in the graph. In step 3, it only needs to check the $M_i$ form, which costs in little time. In step 4, Tarjan's depth-first search algorithm to find the strong components of a directed graph is in linear time of its edges, that is, the bound to the optimal $\textit{O}(N + \Phi)$. In all, the worst-case complexity of the algorithm becomes $\textit{O}(pN\Phi)$. $\Box$
\end{proof}

The $\Sigma$ of sparse pattern $S$ for a given DAEs system can be block triangulated by the triangularization algorithm mentioned as above. Therefore, the $\Sigma$ with block triangular form is as follows:
\begin{equation}\label{MBTF}
\Sigma=
\begin{bmatrix}
  \Sigma_{1,1} & \Sigma_{1,2} & \cdots &\Sigma_{1,p} \\
         & \Sigma_{2,2} & \cdots & \Sigma_{2,p} \\
         &  & \ddots & \vdots \\
         &  &  & \Sigma_{p,p}
\end{bmatrix},
\end{equation}
where the elements in the blanks of $\Sigma$ are $-\infty$. In general, there are three possible cases for $\Sigma$ of sparse pattern $S$ as follows, see Figures 1, 2 and 3.
\begin{figure*}[htb]
\begin{minipage}[t]{0.33\linewidth}
\centering
\includegraphics[width=0.9\textwidth]{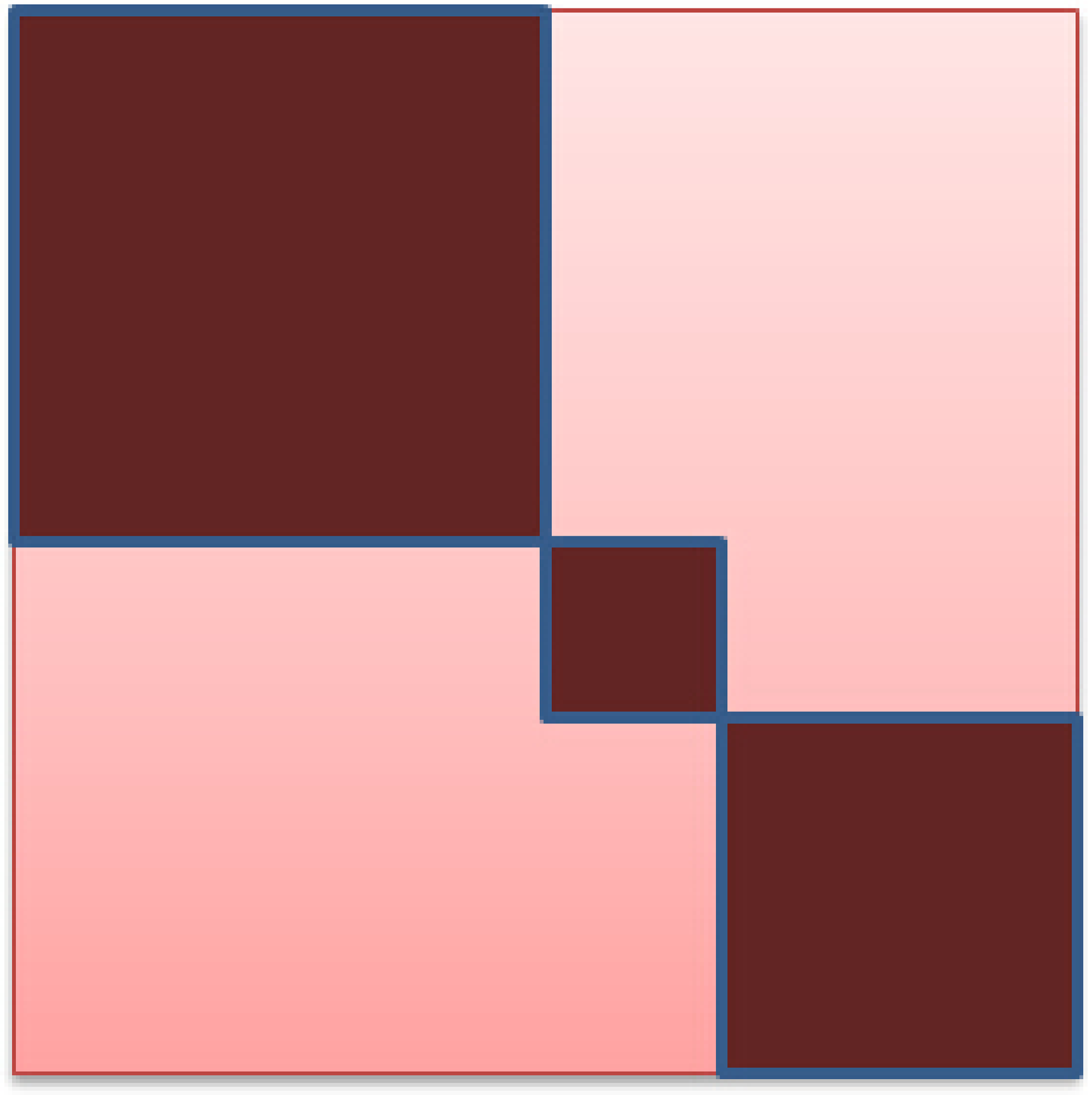}
\caption{irreducible $\Sigma$}
\label{fig:ibtf}
\end{minipage}
\begin{minipage}[t]{0.33\linewidth}
\centering
\includegraphics[width=0.9\textwidth]{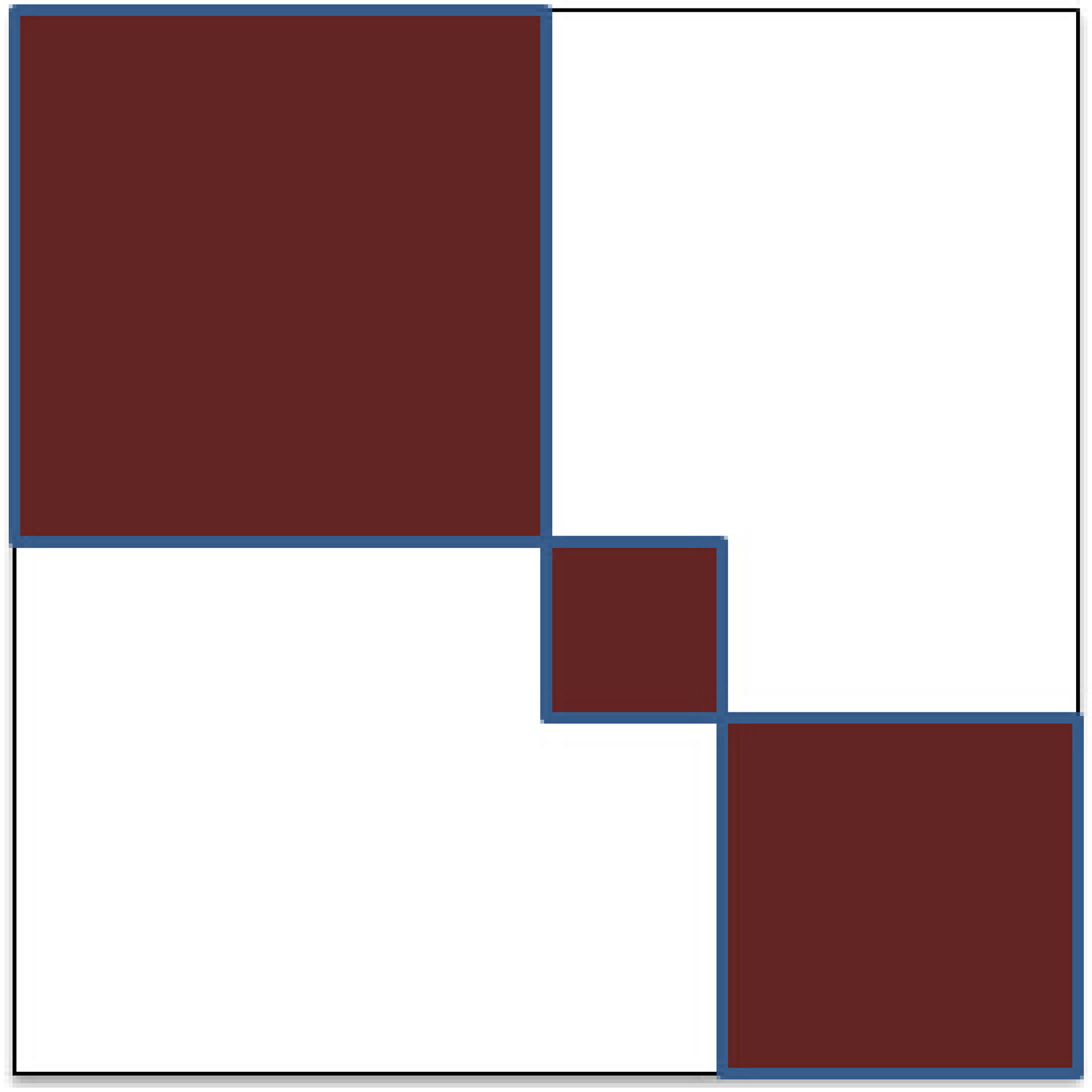}
\caption{$\Sigma$ with only diagonal blocks}
\label{fig:dbtf}
\end{minipage}
\begin{minipage}[t]{0.33\linewidth}
\centering
\includegraphics[width=0.9\textwidth]{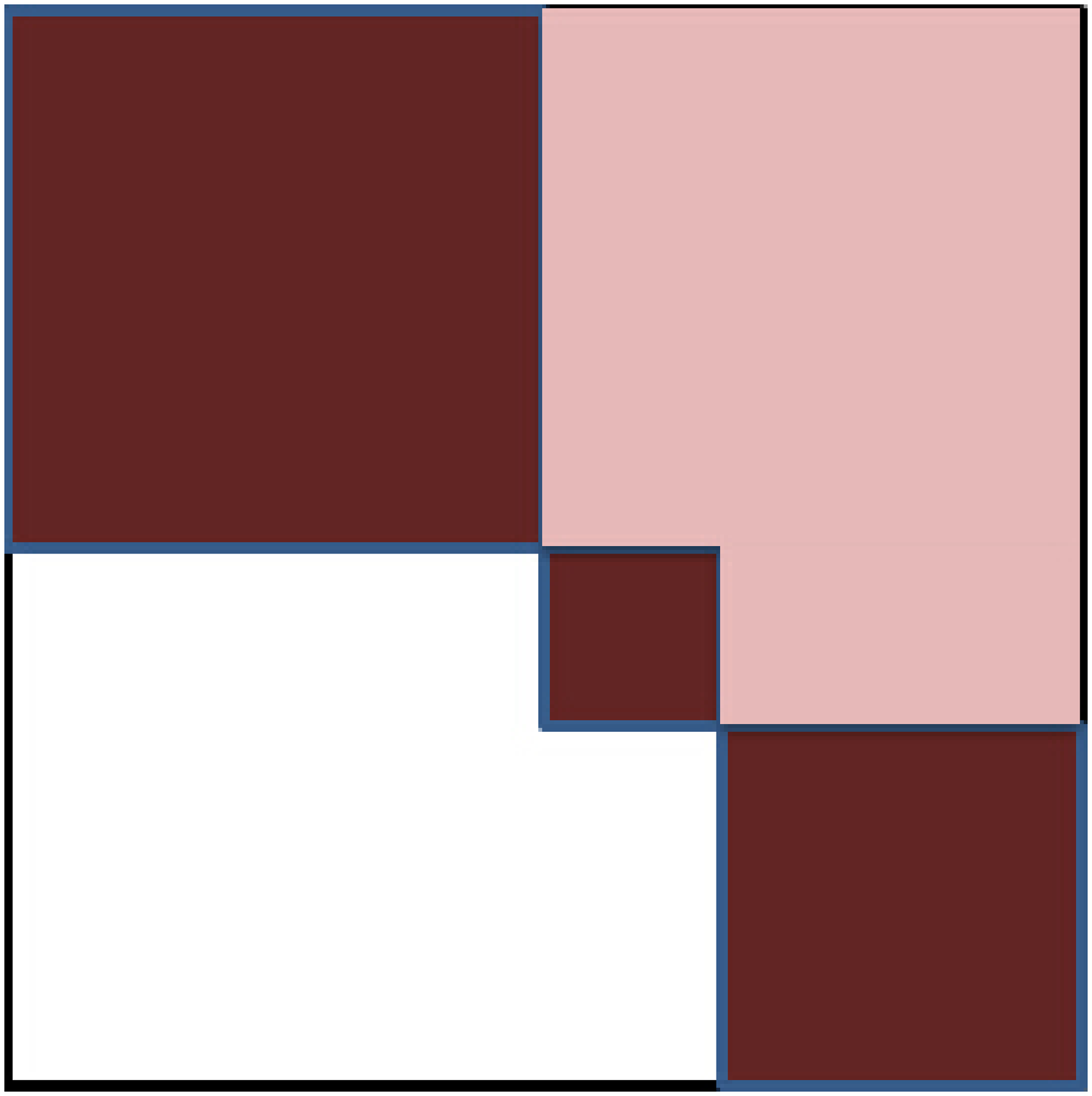}
\caption{$\Sigma$ with block triangular form}
\label{Fig:BTF}
\end{minipage}
\end{figure*}

From Figures $1$, $2$ and $3$, we have the following observations.
\begin{itemize}
 \item For the $\Sigma$ of sparse pattern $S$, there may exist some strong coupling terms with the block triangular form in Figure 1 other than Figures 2, 3. It cannot be permuted to a nontrivial BTF (one with $p > 1$), whose matrix is said to be irreducible, although there exists the maximum number of nonzeros on the diagonal of $\Sigma$. If the $\Sigma$ can be permuted to the BTF (\ref{MBTF}) with $p>1$ in Figures 2, 3 is called to be reducible. Therefore, if the $\Sigma$ contains one transversal, that is, $S$ for the $\Sigma$ contains one transversal, then the $\Sigma_{ij}$ is irreducible for $i, j = 1, 2, \cdots, p$. If the $\Sigma$ for the given DAEs system contains some transversal, the BTF of $\Sigma$ can be found by the triangularization algorithm as above, and its diagonal blocks are square and irreducible.
 \item For Figure 2, each element of each diagonal block is on a transversal from Remark \ref{lem:irreducibility}(b)(iii) in the absence of coupling. It can lead to the parallelization for computing the transversal on each block naturally. The $\Sigma$ contains a remarkable property, which of difference between global offsets of $\Sigma$ and local offsets on each block must be the constant.
 \item For Figure 3, the $\Sigma$ can be permuted to the block upper triangular form. It will be influenced by the top right blocks for computing the transversal, since the lower left blocks are blank. Moreover, it needs the block fixed point iteration method with parameter for computing the offsets of $\Sigma$ from top to bottom in sequence. In particular, for all column $j$ if $\sigma_{ij}$ in the top right blocks are less than the elements of the same columns in the diagonal blocks all the time, then it is similar to Figure 2.
\end{itemize}

\subsection{Shortest augmenting path algorithm}
\label{sec:sapa}
In this subsection, we present the shortest augmenting path algorithm to find a maximum value transversal in the sparse pattern $S$ of signature matrix $\Sigma = (\sigma_{ij})$ from \cite{JV1987}. A transversal of $S$ is an $n$ element subset $T$ of $S$ with just one element in each row and column. The DAEs system is structurally well-posed if there exists a $T$, all of whose $\sigma_{ij}$ are finite, else structurally ill-posed. The problem is to find a maximum value transversal, which makes $||T|| = \sum_{(i,j)\in T} \sigma_{ij}$ as large as possible. Meanwhile, it is equivalent to finding a maximum cardinality matching in a bipartite graph whose incidence matrix is the signature matrix, and can be solved by the Algorithm \ref{alg:MVT}.

\begin{algorithm}[!h]
\caption{(Shortest augmenting paths based algorithms)}\label{alg:MVT}
\begin{algorithmic}[1]
\REQUIRE signature matrix $\Sigma = (\sigma_{ij})$.
 \ENSURE either output a maximum value transversal $T$ for $\Sigma$  or give a nonexistent error.
\end{algorithmic}\label{algo:integer relation}

Step 1: {\em Initiation.}
\begin{algorithmic}[1]
\STATE {\em Column reduction.}
\FOR{$j=n$ to $1$}
\STATE find minimum value over rows
\IF {minimum value in a row appears at an unassigned column}
\STATE initialize the assignment column
\ENDIF
\ENDFOR
\STATE {\em Reduction transfer.}
\FOR{each assigned row $i$}
\STATE transfer from unassigned to assigned rows to a column $k$
\ENDFOR
\STATE {\em Augmenting row reduction.}
\end{algorithmic}
 \begin{algorithmic}[(1)]
\REPEAT
\STATE finding augmenting paths starting in an unassigned row $i$
\IF {$j$ is unassigned}
 \STATE  augmentation solution is from the alternating path
\algstore{myalg}
\end{algorithmic}
\end{algorithm}

\begin{algorithm}[!h]
 \begin{algorithmic}[(1)]
\algrestore{myalg}
 \ELSE
 \STATE reassigning column $j$ to row $i$, reduction is transferred
\ENDIF
\UNTIL{no reduction transfer or augmentation}.
\end{algorithmic}
Step 2: {\em Augmentation.}
\begin{algorithmic}[1]
\REPEAT
\STATE\label{algostep:corner} {\em Augment solution for each free rows.}
\STATE construct the auxiliary network and determine from an unassigned row $i$
to an unassigned column $j$
\STATE find an alternating path by the modified Dijkstra's shortest path method, which is used to augment the solution
\STATE {\em Adjust the signature matrix.}
\STATE update column values, and reset row and column assignments along the alternating path
\UNTIL{no unassigned rows.}
\end{algorithmic}
Step 3: \textbf{return} the sequence of rows, and the corresponding column indices.
\end{algorithm}

\begin{remark}
The complexity of the classical algorithm for finding maximum value transversal can be done at $\textit{O}(n^3)$ operations based on the Hungarian method. For instance, Balinski \cite{B1985} presented a signature methods for the assignment problem, which considers feasible dual solutions corresponding to trees in the bipartite graph of row and column nodes. The shortest augmenting path algorithm can be reduced to $\textit{O}(n^2\log{n})$ operations by using priority queues. By fully exploiting monotonically nondecreasing over all augmentations for the minimum distance, it can be reduced the average time per augmentation to $\textit{O}(n)$. For the sparse case, the whole complexity is more like $\textit{O}(n^2)$ operations. Moreover, we give a necessary preprocessing to the negative of the original $\Sigma$, since the shortest augmenting path algorithm is to compute the minimum-cost network flow problem. A parallel version of Algorithm \ref{alg:MVT} can speed up the MVT finding procedure \cite{BMPT1991}.
\end{remark}

\begin{example}\label{exam1}
We present an index-five model of the planar motion of a crane to illustrate the Algorithm \ref{alg:MVT}. The model is discussed in \cite{C1995}, see Figure \ref{fig:cra}. It is aimed to determine the horizontal velocity $u_1$ of a winch of mass $M_1$, and the angular velocity $u_2$ of the winch so that the attached load $M_2$ moves along a prescribed path $p_1(t)$ and $p_2(t)$.
\begin{figure}[h!]
  \centering
    \includegraphics[width=0.5\textwidth]{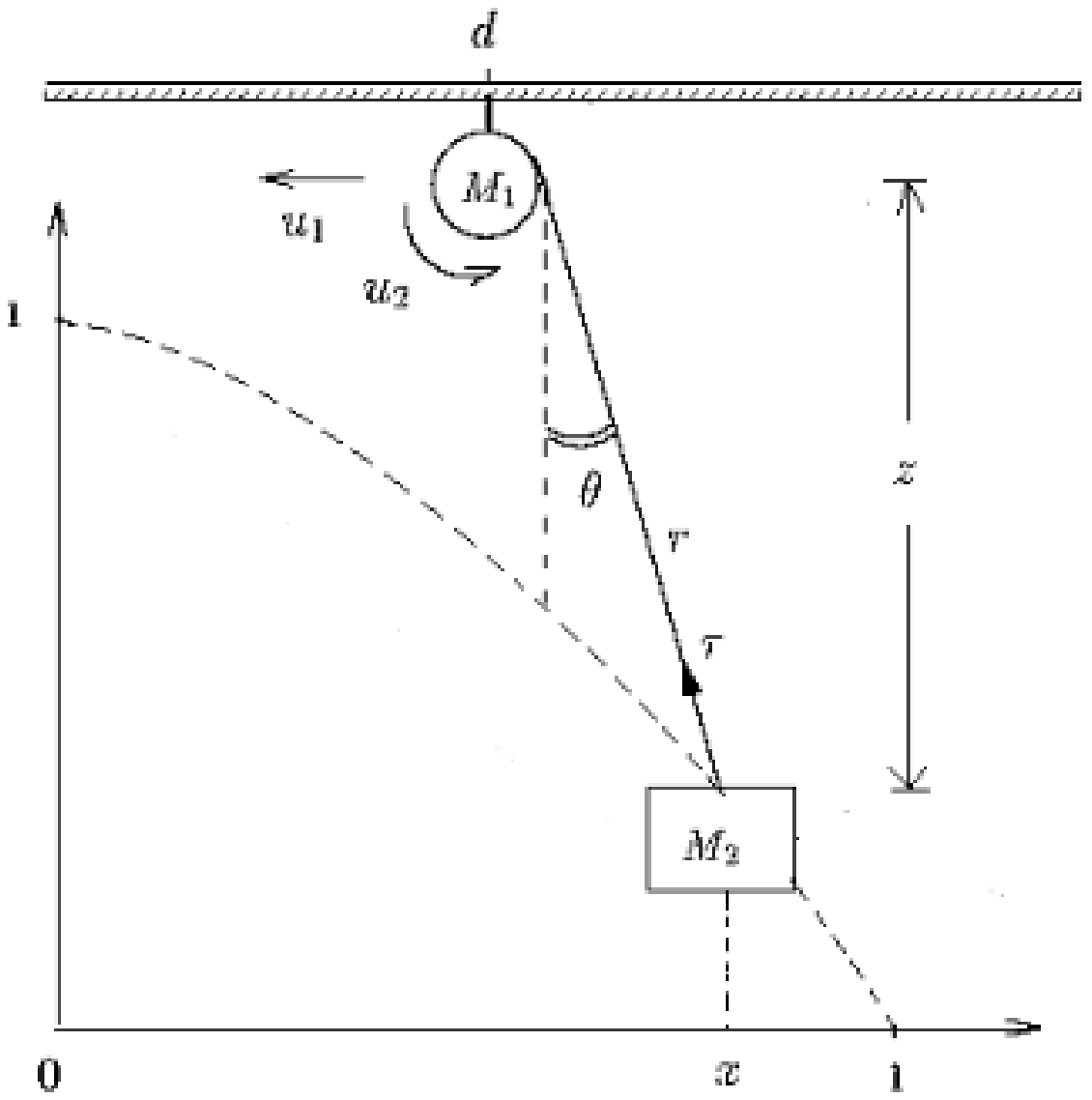}
     \caption{Control of a Crane}
     \label{fig:cra}
\end{figure}

The resulting DAEs are given by the Newton's second law of motion as follows,
\begin{equation}\label{exam1orig}
\left.
\begin{aligned}
0 &= f_1 = M_2\ddot{x}+\tau sin \theta, \ \ \ \ \ \ \ \ \ \ \ \ \ \ \ \ \ \ \ \ \  \ \   0 = f_2 = M_2\ddot{z}+\tau cos \theta - mg, \\
0 &= f_3 = M_1\ddot{d}+ C_1 \dot{d} - u_1 - \tau sin \theta, \ \ \ \    0 = f_4 = J\ddot{r} + C_2\dot{r} + C_3u_2 - C_3^2\tau, \\
0 &= f_5 = r sin\theta +d -x, \ \ \ \ \ \ \ \ \ \ \ \ \ \ \ \ \ \ \ \ \ \  0 = f_6 = r cos \theta - z, \\
0 &= f_7 = x -p_1(t), \ \ \ \ \ \ \ \ \ \ \ \ \ \ \ \ \ \ \ \ \  \ \ \ \ \ \ \ \ 0 = f_8 = z -p_2(t),     
\end{aligned} \ \ 
\right \}
\end{equation}
where $(x, z)$ is the location of the load, $\tau$ is the tension in the cable, $\theta$ is the angle of the cable with the vertical, $J$ is the moment of inertia of the winch, $d$ is the trolley location and $r$ is the cable length, $C_1 > 0$, $C_2 > 0$ $C_3 > 0$, $g > 0$ and $m > 0$ are constants.

From the definition of $\Sigma$ in Section \ref{sec:purp}, labeled by equations and dependent variables, is
\begin{equation}
\label{exam:orgi}
\begin{array}{ll}
\Sigma = \begin{blockarray}{ccccccccc}
~ & x& z & d & r & \theta & \tau & u_1 & u_2 \\
\begin{block}{c[cccccccc]}
f_{1} & 2 & ~ & ~ &  ~  & 0 & 0 & ~ & ~  \\
f_{2} & ~ & 2 & ~ &  ~ & 0  & 0 & ~  & ~  \\
f_{3} & ~ &  ~& 2 &  ~ &  0  & 0 & 0   & ~ \\
f_{4} & ~ & ~ & ~& 2 &  ~ &0 & ~  & 0   \\
f_{5} & 0 & ~ & 0 & 0 & 0 & ~ &  ~  & ~  \\
f_{6} & ~ & 0 & ~ & 0  & 0 & ~&  ~  & ~  \\
f_{7} & 0 & ~ & ~ & ~ & ~ & ~ & ~ & ~  \\
f_{8} & ~ & 0 & ~& ~& ~ & ~ & ~&  ~\\
\end{block}
\end{blockarray} 
\end{array},
\end{equation}
where a blank denotes $-\infty$.

From Section \ref{sec:btf}, we have the block triangularization for (\ref{exam:orgi}) as follows:
\begin{equation*}
\centering
\begin{array}{ll}
\Sigma = \begin{blockarray}{ccccccccc}
~ & u_2 & u_1 & d & r & \tau &\theta & z & x  \\
\begin{block}{c[cccccccc]}
f_{4} &0 & ~ & ~ &  2 \arrayvline  & 0 & ~ & ~ & ~ \\ 
f_{3} &~ &0 & 2 &  ~ \  \arrayvline & 0 & 0 & ~  & ~  \\
f_{5} &~ &  ~& 0 &  0 \arrayvline &  ~  & 0 & ~   & 0 \\
f_{6} &~ & ~ & ~& 0\arrayvline  &  ~ &0 & 0  & ~   \\
\cline{2-7}
f_{1} & ~ & ~ & ~  & ~\ \arrayvline & 0 & 0\arrayvline &  ~  & 2  \\ 
f_{2} & ~ & ~ & ~ & ~ \ \arrayvline & 0 & 0\arrayvline &  2  & ~  \\
\cline{6-9}
f_{8} & ~ & ~ & ~ & ~ & ~ & ~\ \arrayvline  & 0 & ~  \\
f_{7} & ~ & ~ & ~& ~& ~ & ~\ \arrayvline  & ~&  0\\
\end{block}
\end{blockarray} 
\end{array},
\end{equation*}
which is similar to the form of Figure \ref{Fig:BTF}.

From Section \ref{sec:sapa}, we can get the MVT by Algorithm \ref{alg:MVT} : ($f_1$, $\tau$), ($f_2$, $\theta$),
($f_3$, $u_1$), ($f_4$, $u_2$), ($f_5$, $d$), ($f_6$, $r$), ($f_7$, $x$), ($f_8$, $z$).
\end{example}

\section{An extended $\Sigma$-method}
We generalize the signature matrix method to the large scale DAEs system with block fixed point iteration. The basic idea mainly includes the preprocess of the original system with the block triangularization, by adopting the shortest augmenting path algorithm for finding maximum value transversal and solution of dual problem with block fixed point iteration. Without loss of generality, we assume that the DAEs system be structurally well-posed, that is, there exists a transversal. The major steps are as follows.
\begin{description}
  \item[step 1] Form the $n \times n$ signature matrix $\Sigma = (\sigma_{ij})$ of the DAEs system with sparse pattern $S$.

  \item[step 2] Call the block triangularization algorithm in Section 2.2 for dealing with $\Sigma$, and obtain the block form $\Sigma_{ij}$ for $i, j = 1, 2, \cdots, p$.

  \item[step 3] Solve an assignment problem to find the MVT $T_k (k =1, 2, \cdots, p)$ from each block $\Sigma_{kk}$ by the shortest augmenting path algorithm separately.

  \item[step 4] Determine the local canonical offsets of the problem, which are the vectors $\mathbf{(c_k)}_{1\leq k\leq p}= (c_i)_{1\leq i\leq n_k}, \mathbf{(d_k)}_{1\leq k\leq p}=(d_j)_{1\leq j\leq n_k}$ such that $d_j - c_i \geq \sigma_{ij}$, for all $1\leq i, j\leq n_k$, and $d_j - c_i = \sigma_{ij}$ when $(i,j)\in T_k$. The global offsets $\mathbf{(c)}= (c_i)_{1\leq i\leq n}, \mathbf{(d)}=(d_j)_{1\leq j\leq n}$ such that $d_j - c_i \geq \sigma_{ij}$, for all $1\leq i, j\leq n$, and $d_j - c_i = \sigma_{ij}$ when $(i,j)\in T$. It is well known that the difference between global and local offsets is the constant on each block when the irreducible fine BTF is of the Jacobian sparsity pattern \cite{Pryce2014}.

This problem can be formulated as the dual of  (\ref{Lap:Primal}) in the variables $\mathbf{c} = (c_1, c_2, \cdots, c_n$) and $
\mathbf{d} = (d_1, d_2, \cdots, d_n)$, the dual is:
\begin{eqnarray}\label{LPP:Dual}
\begin{array}{ll}
  Minimize & \overline{z} = \sum\limits_j d_j - \sum\limits_i c_i, \label{equ:a}   \\
  subject \ to & d_j - c_i \geq \sigma_{ij} \ for\ all\ (i, j) \in S, \\
  & c_i \geq 0 \ for\ all\ i.
\end{array}
\end{eqnarray}

In order to compute the local canonical offsets of the problem, we apply the fixed point iteration method with parameter to process each irreducible diagonal matrix in block upper triangulated $\Sigma$ from top to bottom in sequence. The theory of block fixed point iteration with parameter method is described in detail by the companion paper \cite{TWQF2014}.

\item[step 5] To verify the success of the index reduction, we need to check whether the $n\times n$ system Jacobian matrix $\mathbf{J}$ is nonsingular, where
\begin{center}
$\mathbf{J}_{ij}=\begin{cases}
 \frac{\partial f_i}{\partial ((d_j-c_i)th\ derivative\ of\ x_j)} \ \ if\ this\ derivative\ is\ present\ in\ f_i,  \\
0 \ \ \ \ \ \ \ \ \ \ \ \ \ \ \ \ \ \ \ \  \ \ \ \ \ \ \ \ \ \ \ \ \
\  otherwise.
\end{cases}$
\end{center}

In this paper the structural index is then defined as:
\begin{center}
$\nu=\max_{i}c_i+\begin{cases}
 0 \ for\ all\ d_j>0 \\
1 \ for\ some\ d_j=0.
\end{cases}$
\end{center}

\item[step 6] Choose a consistent point. If $\mathbf{J}$ is nonsingular at that point, then the
solution can be computed with Taylor series or numerical homotopy continuation techniques in a neighborhood of that point.
\end{description}
It is well known that if $T_k$ is the MVT of $\Sigma_{kk}$ with BTF (\ref{MBTF}), then $T = \bigcup\limits_{k=1}^{p}T_k$ is the MVT of $\Sigma$. In particular, the essential MVT is exactly the union of the diagonal blocks in Figure 2. The parallelization of the case that we have just described can be easily done because it performs the same computations without the necessity of communication between the processors. In Figures 1 and 3, the prior block may influence the posterior blocks for the MVT of $\Sigma$. For the Problem (\ref{LPP:Dual}) of structurally nonsingular DAEs system, the fixed point iteration algorithm can give the same smallest dual-optimal pair with block fixed point iteration algorithm \cite{TWQF2014}.

\begin{theorem}\label{thm:compl}
The above algorithm works correctly as specified and its complexity mainly depends on the block triangularization, finding the MVT and the block fixed point iteration as follows:\\
(a) If $\Sigma$ is irreducible, then the cost of the algorithm is $\textit{O}(\tau + n(nlogn+1)+qn^2)$ operations,\\
(b) If $\Sigma$ is reducible, then the cost of the algorithm is $\textit{O}(p(N\Phi + N^2logN) +qN^2)$ operations,\\
where $q$ is the iteration number of fixed point for computing the global offsets $\mathbf{c}$ in $\Sigma$.
\end{theorem}
\begin{proof}
The Correctness and termination of the algorithm follow from Lemma 3.6 in the companion paper \cite{TWQF2014}.

Its complexity can be easily proved by Theorem 2.4 and Remark 2.5. $\Box$
\end{proof}

\begin{remark}
In steps 3 and 4, for the existing MVT in $\Sigma$, we give the time complexity $\textit{O}(pN^2logN + qN^2)$, which is superior to $\textit{O}(pN^3 + qN^2)$ in \cite{TWQF2014}. In \cite{DZC2008, Pantelides1988}, the whole time complexity of their algorithms is $\textit{O}(n^2 + n\tau + r(n+\tau))$, where $r$ is the number of adjusting the maximum weighted value from negative to zero to find a perfect matching in the bipartite graph. In particular, their index reduction techniques can be viewed as the behavior of signature matrix method by the bipartite graph. However, they can not easily to achieve the block triangularization. Furthermore, Pryce only declared the Algorithm 3.1 in finitely many iterations in \cite{Pryce2001}.
\end{remark}

\begin{example} \label{exam2}
Continue from Example \ref{exam1}, we can mark the MVT with '*' and easily compute its gobal offsets by our extended $\Sigma$-method as follows:
\begin{equation*}
\begin{array}{ll}
\Sigma = \begin{blockarray}{cccccccccccccc}
~       & u_2 & u_1 & d & r & \tau &\theta & z & x  & \mathbf{c}^{(1)}& \mathbf{c}^{(2)}& \mathbf{c}^{(3)}& \mathbf{c}^{(4)}   & \mathbf{c}^{(5)}    \\
\begin{block}{c[cccccccc]ccccc}
f_{4} & 0^* & ~ & ~  &\  2\ \arrayvline  & 0  & ~ & ~ & ~ &0&0& 0& 0& 0\\
f_{3} & ~  &0^* & 2 & \ ~ \ \ \arrayvline & 0 & 0 & ~  & ~  &0&0&0& 0& 0\\
f_{5} & ~ &  ~ & 0^*& \ 0\ \arrayvline&  ~  & 0 & ~   & 0&2&2&2 & 2& 2\\
f_{6} & ~ & ~ & ~ & 0^*  \arrayvline&  ~ & 0 & 0  & ~   &2&2&2& 2& 2\\
\cline{2-7}
f_{1} & ~ & ~ & ~  & ~\ \  \arrayvline& 0^* &\ 0 \  \arrayvline&  ~  & 2  &0&0&\fbox{2}& 2& 2\\
f_{2} & ~ & ~ & ~ & ~\ \ \arrayvline  & 0 & 0^*\arrayvline &  2  & ~  &0&\fbox{2}&2& 2& 2\\
\cline{6-9}
f_{8} & ~ & ~ & ~ & ~ & ~ & ~\ \  \arrayvline & 0^*  & ~  &2&2&\fbox{4}& 4& 4\\
f_{7} & ~ & ~ & ~& ~& ~ & ~ \ \ \arrayvline & ~&  0^* &2&2&2& \fbox{4}& 4\\
\end{block}
\mathbf{d} ^{(1)}& 0& 0 & 2 & 2 &  0 & 0 &2 &2 \\
 \mathbf{d}^{(2)}& 0& 0 & 2 &  2 & 0 &\fbox{2} & 2& 2\\
 \mathbf{d}^{(3)}& 0& 0 & 2 &  2 &\fbox{2} &2 & \fbox{4}& 2\\
 \mathbf{d}^{(4)}& 0& 0 & 2 &  2 & 2 &2 & 4&\fbox{4}\\
 \mathbf{d}^{(5)}& 0& 0 & 2 &  2 & 2 &2 & 4& 4\\
\end{blockarray}
\end{array}.
\end{equation*}
We can get the offsets $\mathbf{c} = (2, 2, 0, 0, 2, 2, 4, 4)$ and
$\mathbf{d} = (4, 4, 2, 2, 2, 2, 0, 0)$, for the variables in the original order given (\ref{exam:orgi}). From step 5, the system Jacobian matrix $\mathbf{J}$  is
\begin{equation}
\label{exam:jac}
\begin{array}{ll}
\mathbf{J}  = \begin{blockarray}{cccccccc}
\begin{block}{[cccccccc]}
 M_2 & 0 & 0 & 0&  \tau cos(\theta)  & sin(\theta) & 0 & 0\\
 0 & M_2 & 0 & 0 & -\tau sin(\theta)  & cos(\theta) & 0  & 0  \\
0 &  0& M_1 &  0 &  0  & 0 & -1   & 0 \\
 0 & 0 & 0& J &  0 &0 & 0  & C_3   \\
 0 & 0 & 1 & sin(\theta) & rcos(\theta) & 0 &  0  & 0  \\
0 & 0 & 0 & cos(\theta)  &-rsin(\theta)  & 0&  0  & 0  \\
 1	 & 0& 0 &  0 & 0 & 0 & 0 & 0  \\
 0 & 1 & 0 & 0 & 0 & 0 & 0 &  0\\
\end{block}
\end{blockarray} 
\end{array},
\end{equation}
Since $det \ \mathbf{J} = cos(\theta)(-\tau cos(\theta)^2C_3 - \tau sin(\theta)^2 C_3)$, whether this system is solvable or not depends on the angle $\theta$, which is easy to check. Consequently, the extended $\Sigma$-method succeeds.
\end{example}

\section{Experimental results}
\begin{figure}[h!]
  \centering
    \includegraphics[width=1\textwidth]{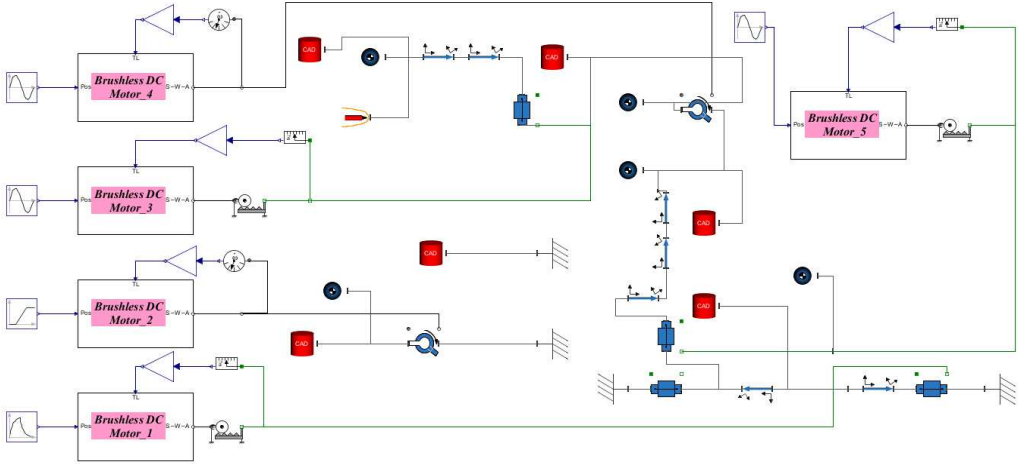}
    \caption{Five-axis Linkage CNC Machine YHV6025 models}
    \label{fig:3DOFM}
\end{figure}

An efficient actual implementation of our algorithm is in
\emph{Maple}. The following examples run in the platform of \emph{Window}
and AMD Athlon(tm) II X4-645 CPU 3.10GHz, 4.00GB RAM. First, we present an industrial application to illustrate our algorithm for the Five-axis Linkage CNC Machine YHV6025 models, see Figure \ref{fig:3DOFM}.

\begin{figure*}[htb]
\begin{minipage}[t]{0.5\linewidth}
\centering
\includegraphics[width=0.9\textwidth]{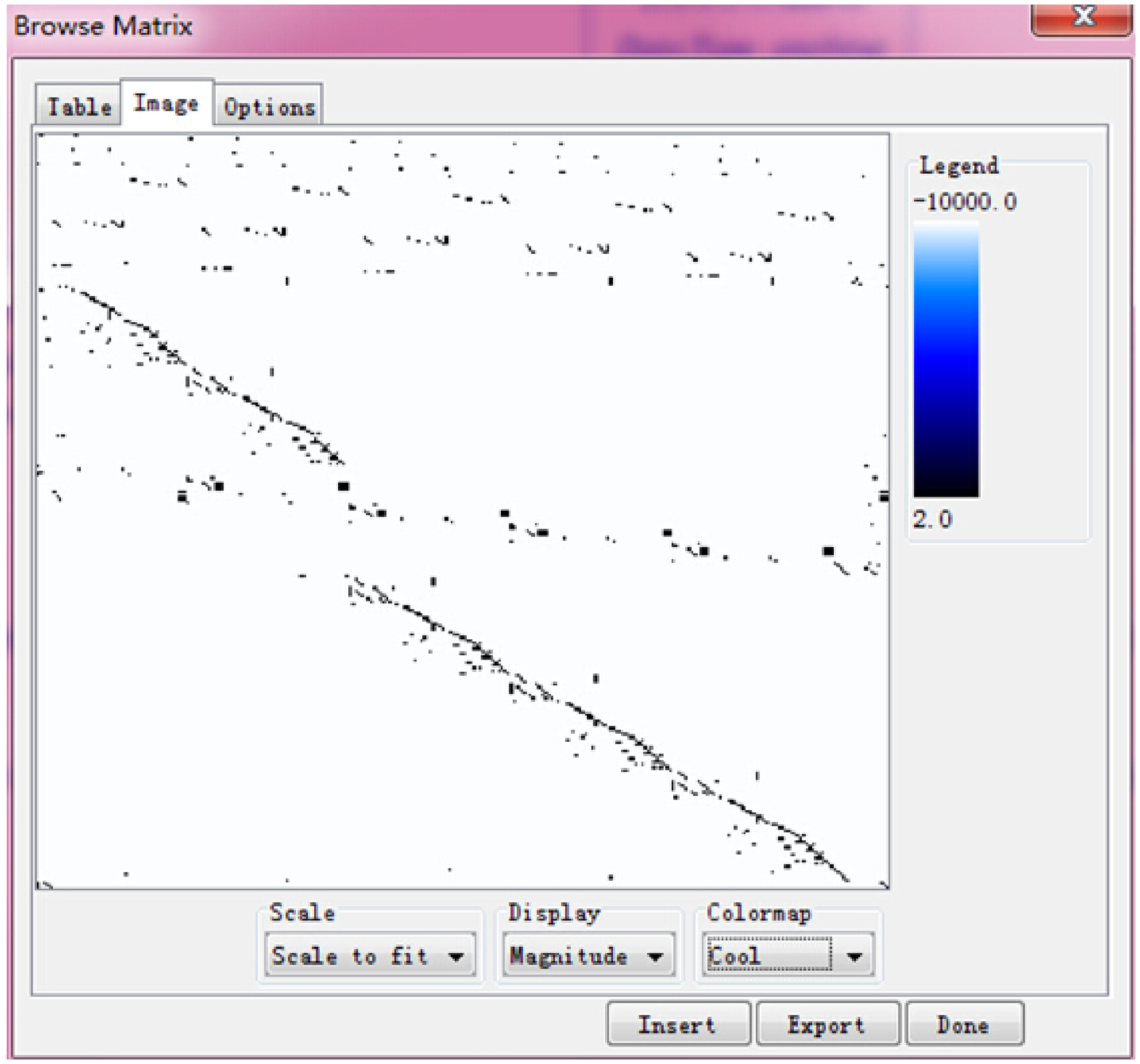}
\caption{the original $\Sigma$}
\label{fig:orig}
\end{minipage}
\begin{minipage}[t]{0.5\linewidth}
\centering
\includegraphics[width=0.9\textwidth]{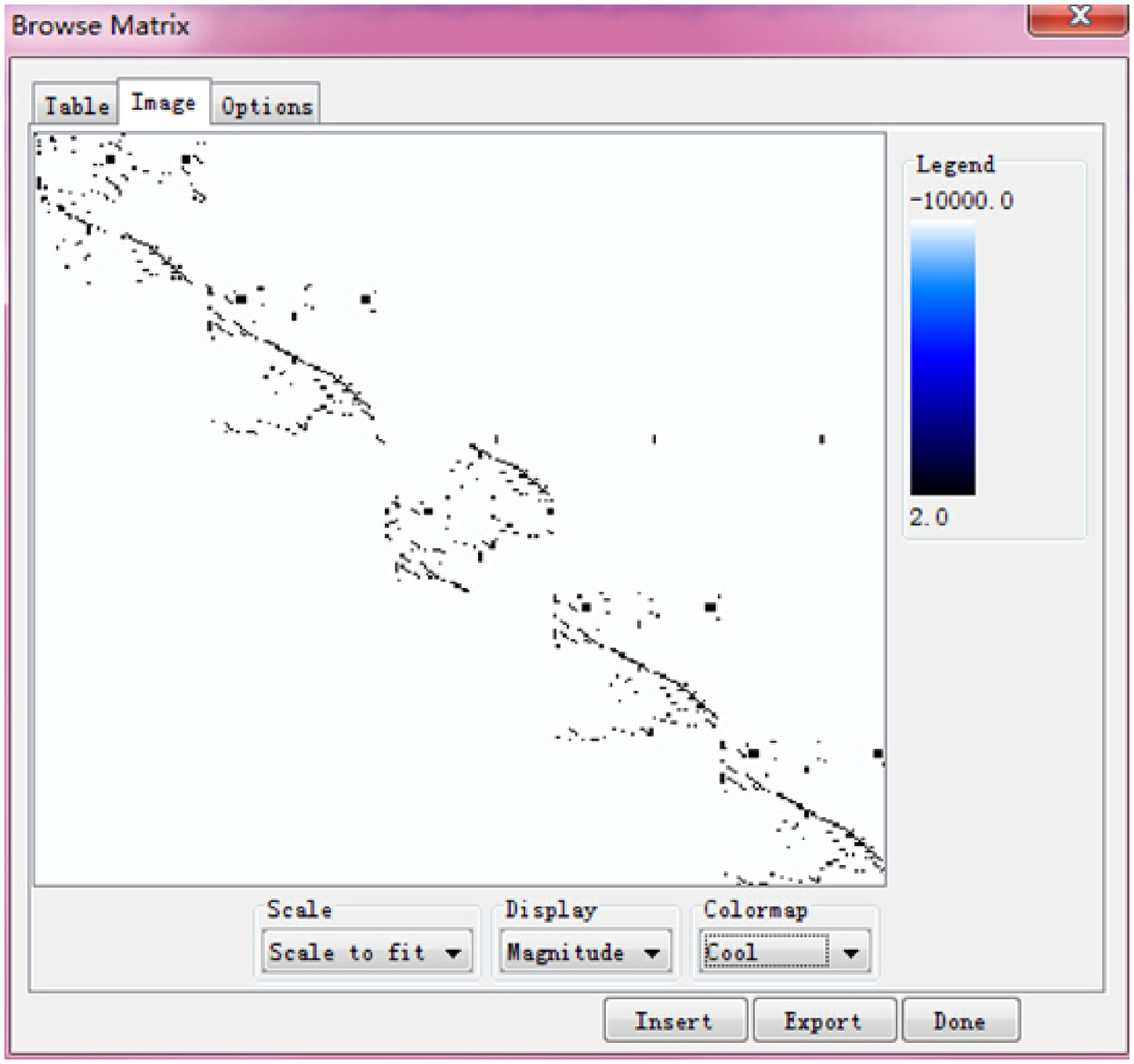}
\caption{the triangulated $\Sigma$}
\label{fig:tria}
\end{minipage}
\end{figure*}

There are 2446 differential and algebraic equations in the original physical system based on the Modelica language tool. We can get 245 differential and algebraic equations through the efficient elimination of trivial equations by means of symbolic simplification. The original sytem $\Sigma$ represents in Figure \ref{fig:orig}. In Figure \ref{fig:tria}, we can obtain the block structure signature matrix for the corresponding $\Sigma$ based on our block triangularization algorithm.

The triangulated signature matrix includes three blocks as follows.
\begin{eqnarray*}
\Sigma =
\begin{bmatrix}
\Sigma^{(1)} \\
& \Sigma^{(2)} & \\
& & \Sigma^{(3)}
\end{bmatrix},
\Sigma^{(1)} =
\begin{bmatrix}
0& \times \\
& \Sigma_{22}^{(1)}
\end{bmatrix},
\Sigma^{(2)}  =
\begin{bmatrix}
0& \times \\
& \Sigma_{22}^{(2)}
\end{bmatrix},
\Sigma^{(3)} =
\begin{bmatrix}
0&\times&\times&\times&\times&\times \\
&0 &\times&\times&\times&\times \\
& &0&\times&\times&\times\\
& & & \Sigma_{44}^{(3)}&\times&\times\\
& & & &\Sigma_{55}^{(3)} &\times\\
& & & & &\Sigma_{66}^{(3)}
\end{bmatrix},
\end{eqnarray*}
where the blank in $\Sigma$, $\Sigma^{(1)}$, $\Sigma^{(2)}$ and $\Sigma^{(3)}$ means -$\infty$, and the orders of signature matrices $\Sigma_{22}^{(1)}$, $\Sigma_{22}^{(2)}$, $\Sigma_{44}^{(3)}$, $\Sigma_{55}^{(3)}$ and $\Sigma_{66}^{(3)}$ are 54, respectively. For each $\Sigma^{(i)}(i=1, 2, 3)$, we present the time for finding the global offsets $\mathbf{(c)}$ and $\mathbf{(d)}$ by the extended signature matrix method (ESMM), and then converting DAEs to the corresponding ODEs and computing the consistent initial values. In Table 1, we compare our method with signature matrix method (SMM) \cite{Pryce2001} and weighted bipartite graph based on index reduction (WBGIR) \cite{DZC2008}.

\begin{table}[ht]
\caption{Time for structural analysis of YHV6025 models} 
\centering
\begin{tabular}{c|c|c|c}
\hline Algorithm & $\Sigma^{(1)}$ &$\Sigma^{(2)}$&$\Sigma^{(3)}$ \\
\hline 
ESMM& 0.06036s&0.05481s&0.13926s\\
SMM& 0.06942s&0.06417s&0.44486s\\
WBGIR& 0.06823s&0.06400s&0.43745s\\[1ex]
\hline 
\end{tabular}
\label{table:nonlin} 
\end{table}
From Table 1, our algorithm is more than three times faster than the other two algorithms for the large scale $\Sigma^{(3)}$. For the small and middle scale $\Sigma^{(1)}$ and $\Sigma^{(2)}$, the advantage of our algorithm is not obvious compared with the other two algorithms. It shows that our algorithm is efficient for the large scale system of DAEs. 

\begin{figure}[!ht]
\begin{subfigure}{.33\textwidth}
  \centering
  \includegraphics[width=1\linewidth]{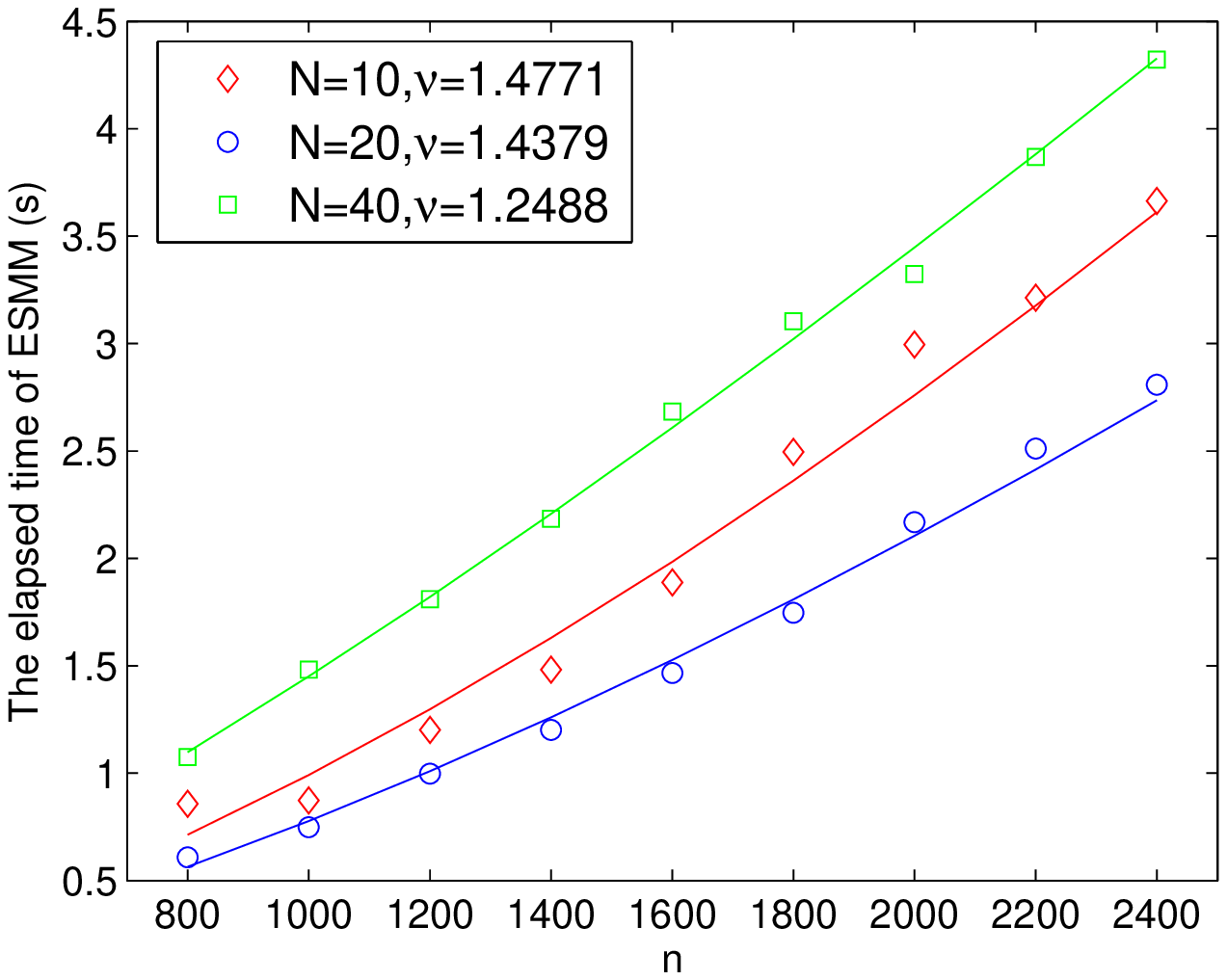}
  \caption{ESMM}
  \label{fig:ESMM}
\end{subfigure}%
\begin{subfigure}{.33\textwidth}
  \centering
  \includegraphics[width=1\linewidth]{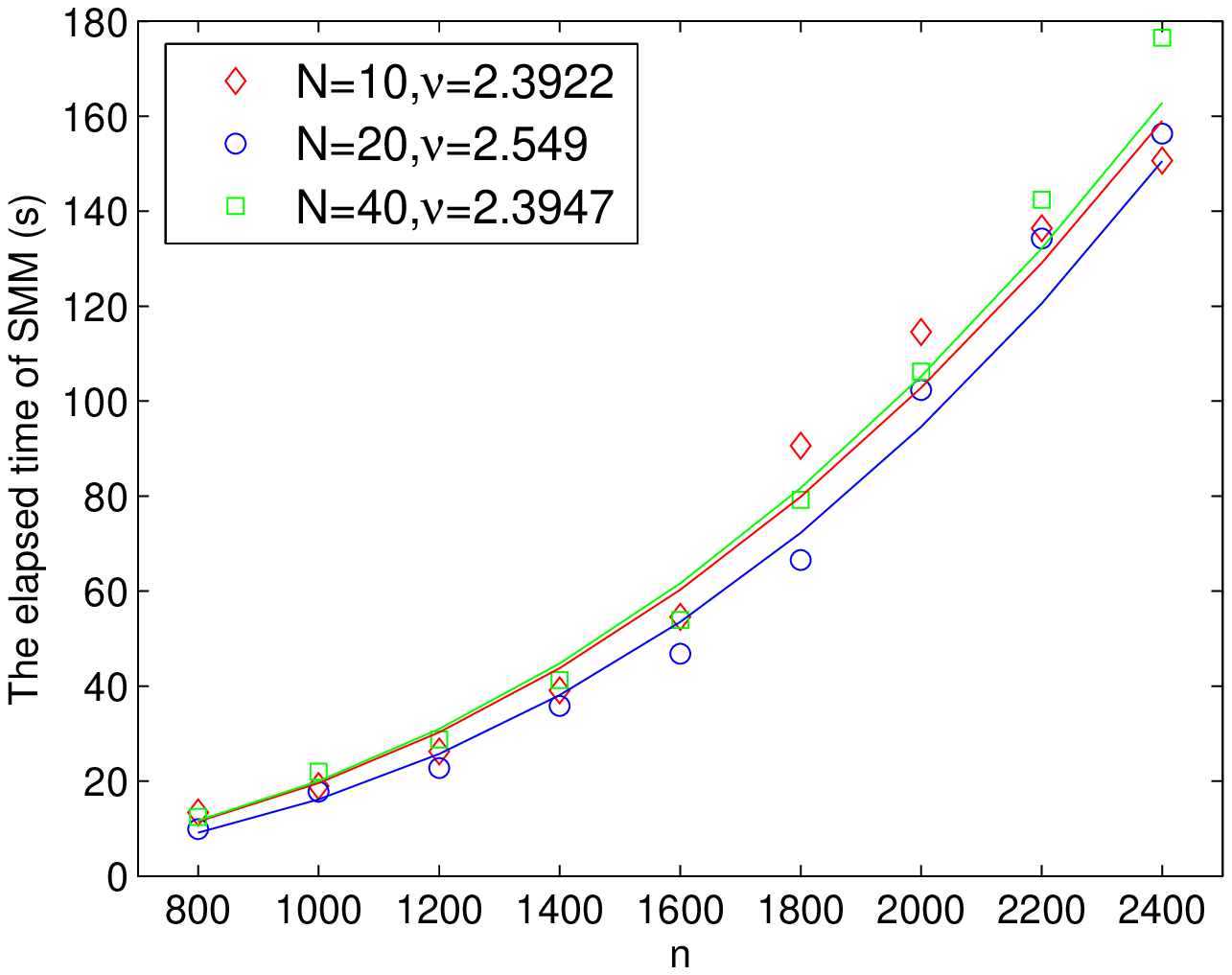}
  \caption{SMM}
  \label{fig:SMM}
\end{subfigure}
\centering
\begin{subfigure}{.33\textwidth}
  \centering
  \includegraphics[width=1\linewidth]{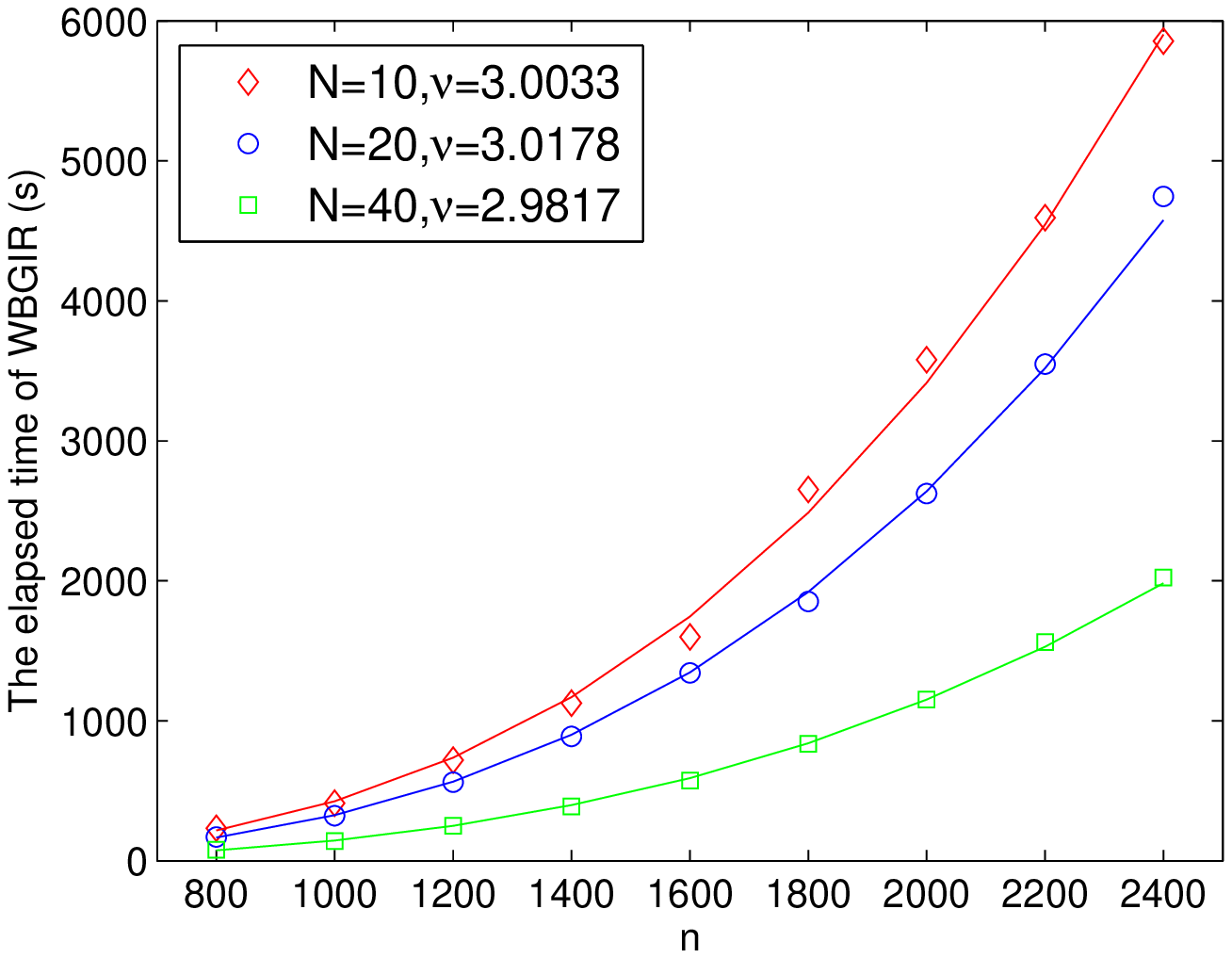}
  \caption{WBGIR}
  \label{fig:WBGIR}
\end{subfigure}
\caption{The elapsed times of structural analysis methods versus $n$. The constants in the legends are the fitted values $\mu$ in $\mu n^\nu$ for different $N$.}
\label{fig:SA}
\end{figure}

\begin{figure}[H]
\begin{subfigure}{.5\textwidth}
  \centering
  \includegraphics[width=1\linewidth]{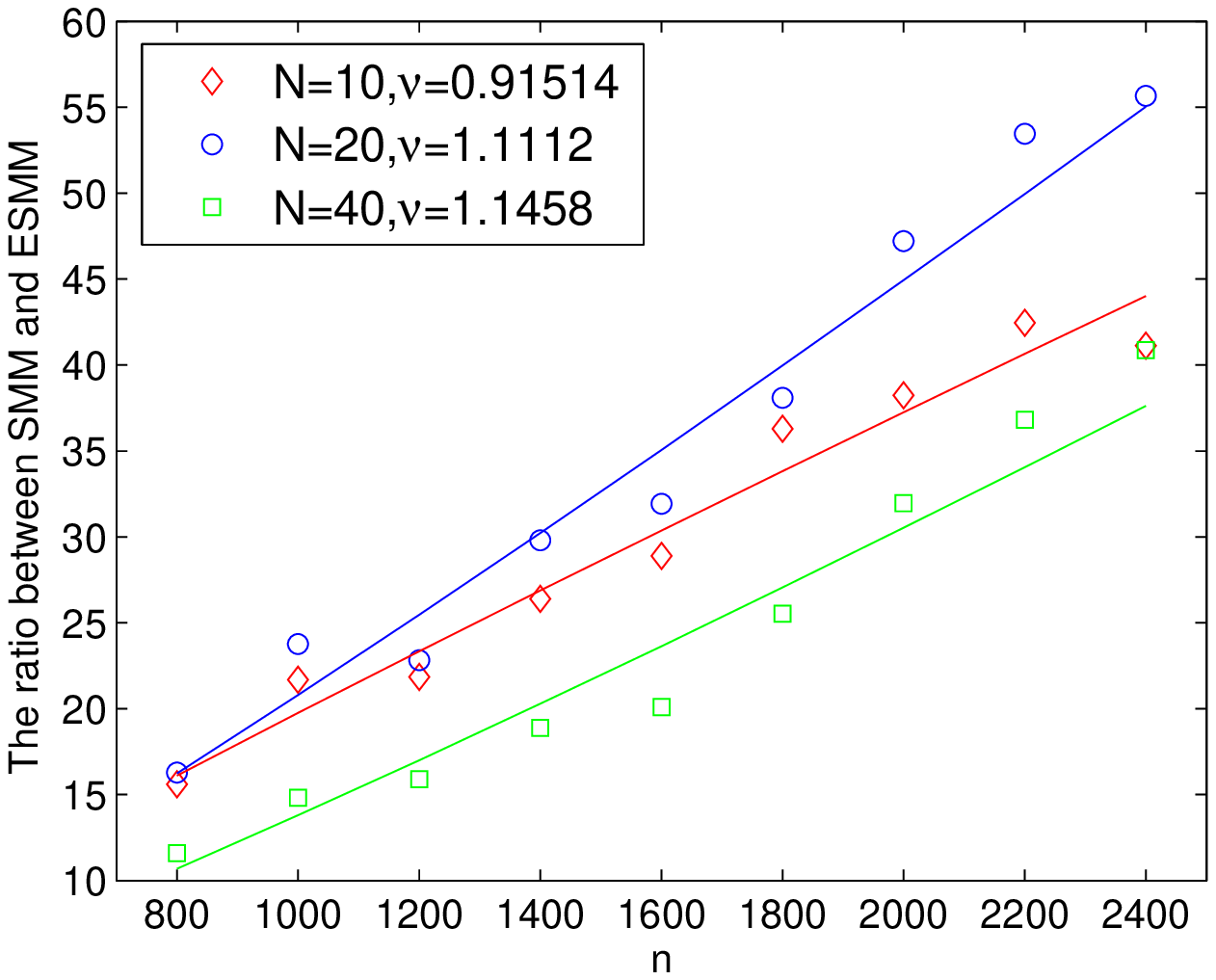}
  \caption{The ratio of elapsed times, i.e., $\frac{\text{SMM's\ elapsed\ time}}{\text{ESMM's\ elapsed\ time}}$, versus $n$}
  \label{fig:Ratio1}
\end{subfigure}%
\begin{subfigure}{.5\textwidth}
  \centering
  \includegraphics[width=1\linewidth]{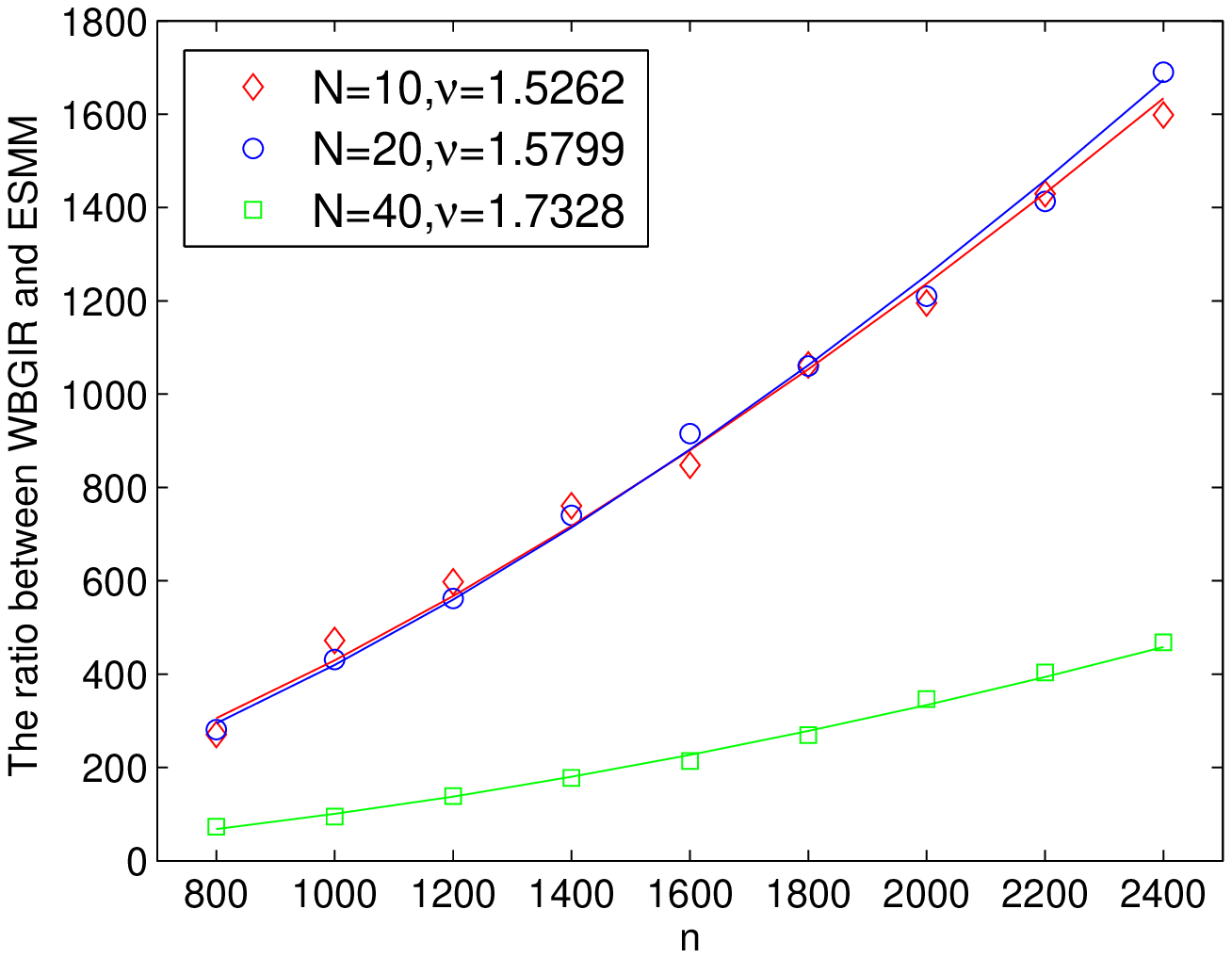}
  \caption{The ratio of elapsed times, i.e., $\frac{\text{WBGIR's\ elapsed\ time}}{\text{ESMM's\ elapsed\ time}}$, versus $n$}
  \label{fig:sfig2}
\end{subfigure}
\caption{Some ratios of elapsed times versus $n$. The constants in the legends are the fitted values $\mu$ in $\mu n^\nu$ for different $N$.}
\label{fig:Ratio}
\end{figure}

We also give some results for large scale systems of differential algebraic equations via numerical simulation experiments.
Assume that the original signature matrices $\Sigma$ of systems have been preprocessed for the BTFs in Section \ref{sec:btf} and
that the order of each diagonal block $\Sigma_{ii}$ in the $n\times n$ matrix $\Sigma$ is $N$ in our experiments.
In general, it is noted that the entries of $\Sigma$ are the differential orders in DAEs. With loss of generality, we can randomly generate the sparse $\Sigma$ for systems
with BTFs properly as follows:
\begin{itemize}
  \item for each element in $\Sigma_{1,1}$, randomly select an integer
   from \{0, 1, 2, 3\} according to its probability from \{0.7, 0.1, 0.1, 0.1\}, respectively; and $\Sigma_{i,i}=\Sigma_{1,1}, i=2,3,\ldots,p$;
  \item for each element in $\Sigma_{1,2}$, randomly select an integer
   from \{-1000, 0, 1, 2\} according to its probability from \{0.925, 0.025, 0.025, 0.025\}, respectively; and $\Sigma_{i,i+1}=\Sigma_{1,2}, i=2,3,\ldots,p-1$;
  \item the rest elements in $\Sigma$ are $-1000$ (means $-\infty$).
\end{itemize}
We test the corresponding random trials for ESMM, SMM and WBGIR
with $N=\{10,20,40\}$ and $n=800:200:2400$, and then calculate their
constants in $\mu \cdot n^{\nu}$ using the standard least-square method.
The elapsed times of three algorithms are shown in Figure \ref{fig:SA}, respectively;
some ratios of elapsed times are given in Figure \ref{fig:Ratio}.

From Figure \ref{fig:SA}, we empirically know that the elapsed times of
ESMM for the large scale nonlinear DAEs with different $N$ are between $O(n)$ and $O(n^{1.5})$, and SMM seems to be $O(n^{2.5})$. However, WBGIR seems to be $O(n^3)$. It is also noteworthy that WBGIR is very time consuming because of its recursive operations.
Figure \ref{fig:Ratio} shows that ESMM can reduce its elapsed time of sparse systems for fixed $N$ by nearly $O(n)$ (i.e., $O(p)$) compared to SMM
and by $O(n^{1.5})$ at least compared to WBGIR. In particular, the experimental results of ESMM are consistent with our complexity analysis in Theorem \ref{thm:compl}. 

\begin{remark}
For the large scale systems of sparse DAEs, we only consider the general reducible case for structural analysis. Moreover, in order to compare these three algorithms fairly, we construct the spare $\Sigma$ of systems with BTFs properly. In practical applications, the signature matrix of sparse pattern $S$ available from the large scale problems are often reducible.
\end{remark}

\section{Conclusion}
\label{}
In this paper, we propose an effective method to compute the large scale DAEs system in practice. We generalize the $\Sigma$-method with BTFs and combine with the block fixed point iteration for the general reducible DAEs systems. In particular, we exploit the shortest augmenting path algorithm for finding maximum value transversal, which can reduce the cost of computation. However, we apply the $\Sigma$-method as a basic tool, which relies heavily on square and sparsity structure. Thus it may be confronted with the drawback that can succeed yet not correctly in some DAEs arising from the specific types \cite{Pantelides1988, Reissig2000}. Some researchers also proposed some new computational techniques for high index nonlinear DAEs solving, such as the direct numerical computation\cite{EH2009, L2014}, index reduction \cite{MS1993, QYFBF2015, TI2008}, and so on. 

In the future, we would like to consider the further research of hybrid symbolic and numerical computation with BTFs for large scale nonlinear systems of DAEs solving in practical industrial applicaitons, such as numerical algebraic geometry. 
\section*{Acknowledgement}

This research was partly supported by China 973 Project NKBRPC-2011CB302402, the National Natural Science Foundation of
China (No. 61402537, 91118001), the West Light Foundation of the Chinese Academy of Sciences, the China Postdoctoral Science
Foundation funded project (No. 2012M521692), and the Open Project of Chongqing Key Laboratory of Automated Reasoning and Cognition (No. CARC2014004).

\end{document}